\documentclass[reqno]{amsart}
\usepackage{amssymb,amsmath}
\usepackage[english]{babel}
\usepackage{amsfonts}
\usepackage{color}
\usepackage{amsthm}
\usepackage[colorlinks=true,linkcolor=NavyBlue, citecolor=NavyBlue,urlcolor=NavyBlue]{hyperref}
\usepackage{graphicx}
\usepackage{mathtools}
\usepackage[usenames,dvipsnames]{pstricks}
\usepackage{bm}
\usepackage{amsmath}
\usepackage{soul}
\usepackage{cancel}

\usepackage{xypic}
\usepackage[draft]{fixme}
\DeclareMathOperator{\AGL}{AGL}

\DeclareMathOperator{\Sym}{Sym}

\DeclareMathOperator{\N}{\mathbb{N}}
\DeclareMathOperator{\mex}{mex}

\newcommand\deq{\mathrel{\stackrel{\makebox[0pt]{\mbox{\normalfont\tiny def}}}{=}}}

\newcommand{\Size}[1]{\left\lvert #1 \right\rvert}
\newcommand{\Span}[1]{\left\langle\,#1\,\right\rangle}
\newcommand{\Set}[1]{\left\lbrace #1 \right\rbrace}

\newcommand{\puncture}[2]{\mathop{\vee}\!\left[#1;#2 \right]}

\let\phi\varphi
\theoremstyle{plain}
\newtheorem{theorem}{Theorem}%[section]

\newtheorem{proposition}[theorem]{Proposition}
\newtheorem{corollary}[theorem]{Corollary}

\theoremstyle{remark}
\newtheorem{remark}{Remark}%[theorem]

\theoremstyle{definition}
\newtheorem{definition}[theorem]{Definition}
\newtheorem*{definition*}{Definition}
\newtheorem{example}[theorem]{Example}
\newtheorem*{notation*}{Notation}
\def\hl#1{{\color{NavyBlue}#1}}
\def\hll#1{{\textbf{#1}}}
\numberwithin{equation}{section}

\title{Unrefinable partitions into distinct parts in a normalizer chain}

{

\author[R.~Aragona]{Riccardo Aragona} 
\author[R.~Civino]{Roberto Civino}
\author[N.~Gavioli]{Norberto Gavioli}
\author[C.~M.~Scoppola]{Carlo Maria Scoppola}
\address[R.~Aragona, R.~Civino, N.~Gavioli]{
	Dipartimento di Ingegneria e Scienze dell'Informazione e Matematica \\ 
	Universit\`a degli Studi dell'Aquila\\
	via Vetoio\\
	I-67100 Coppito (AQ)\\
	Italy}
\email[Riccardo Aragona]{riccardo.aragona@univaq.it}
\email[Roberto Civino]{roberto.civino@univaq.it}
\email[Norberto Gavioli]{norberto.gavioli@univaq.it}
\email[Carlo Maria Scoppola]{carlo.scoppola@univaq.it}

\thanks{All the authors are members of INdAM-GNSAGA (Italy) and of the
  group  ``Crittografia   e  Codici''   of  the   Italian  Mathematical
  Union. R.  Civino is  partially funded by  the Centre  of Excellence
  EX-EMERGE at University of L'Aquila.  N. Gavioli is a member of
  the Centre of Excellence EX-EMERGE at University of L'Aquila.  }


\begin{document}

\begin{abstract}
In a recent paper on a study of the Sylow 2-subgroups of the symmetric 
group with $2^n$ elements it has been show that the growth of the first $(n-2)$ 
consecutive indices of a certain 
normalizer chain is linked to the sequence of 
partitions of integers into distinct parts. Unrefinable partitions into distinct parts are those
 in which no part $x$ can be replaced with integers
 whose sum is $x$ 
 obtaining a new partition into distinct parts. 
We prove here that the $(n-1)$-th index of the previously mentioned chain is related to the number 
of unrefinable partitions into distinct parts satisfying a condition on the minimal excludant.

\end{abstract}

\maketitle

\section{Introduction}
The sequence $(b_j)$ of the number of partitions of integers into distinct parts has been 
extensively studied
in past and recent years and is a well-understood integer sequence~\cite{euler1748introductio,Andrews1994} appearing in different 
areas of mathematics. 
Recently,  triggered by a problem in  algebraic cryptography related to translation subgroups in the      symmetric     group      with      $2^n$
elements~\cite{cds06,Aragona2019,Civino2019,calderini2017translation},
it has been shown that such a sequence is related to the growth of
the indices of consecutive terms in a chain of normalizers~\cite{aragona2021rigid}.
More precisely, let $\Sigma_n$ be  a  Sylow
$2$-subgroup  of  the symmetric group $\Sym(2^n)$ and $T$ be an elementary abelian regular subgroup of $\Sigma_n$. Defining
$N_n^0 = N_{\Sigma_n}(T)$ and recursively $N_n^i = N_{\Sigma_n}(N_{n}^{i-1})$, the authors proved that 
  the   number
\[\log_{2}\Size{N^{i}_n  : N^{i-1}_n}\]  is independent  of $n$  for
$1\le  i\le n-2$,  and is  equal  to the  $(i+2)$-th  term of  the
sequence $(a_j)$ of the partial sums  of  $(b_j)$ (cfr. Table~\ref{tab:two}).

\noindent The result is obtained
by proving that  the terms in the chain are \emph{saturated} subgroups, i.e.\ generated by \emph{rigid commutators},
a family of left-normed commutators involving  a special set of generators of $\Sigma_n$.
We invite the interested reader to refer to Aragona et al.~\cite{aragona2021rigid}, where rigid commutators  and saturated groups are introduced and described in detail.

When $i > n-2$, the behavior of the chain does not seem to show any recognizable pattern and the study of its combinatorial nature is still open. 
Nonetheless, further investigations on experimental evidences led us to notice that the first exception to the rule, i.e.\
$\log_{2}\Size{N^{n-1}_n  : N^{n-2}_n}$, is linked to the number of partitions of $n$ into distinct parts which do not admit a \emph{refinement}, 
i.e.\ those partitions $n=p_1+p_2+\ldots+p_t$ where no integer $p_i$ can be replaced
by $l \geq 2$ distinct integers $q_1,q_2,\ldots,q_l$, whose sum 
is $a_i$ and such that the resulting partition 
\[n=p_1+p_2+\ldots+p_{i-1}+(q_1+q_2+\ldots q_l)+p_{i+1}+\ldots+p_t\] is still a partition into distinct parts.
In our knowledge, \emph{unrefinable} partitions have not been investigated so far. 

We prove here (cfr. Theorem~\ref{prop:elementsofthen-1normalizer}) that a transversal of $N_{n}^{n-2}$ in $N_n^{n-1}$ is made of rigid commutators which are 
 in one-to-one correspondence with unrefinable partitions 
into distinct parts with a condition on their minimal excludants, where the minimal excludant of a partition is the least positive integer which does not appear in the partition. \\

\begin{table}[]
		\centering
	{\renewcommand{\arraystretch}{1.3}
		\begin{tabular}{c||c|c|c|c|c|c|c|c|c|c|c|c|c|c|c}
			$j$ & $0$ & $1$ & $2$ & $3$ & $4$ & $5$ & $6$ & $7$ & $8$ & 9& 10 &11&12&13&14\\
			\hline\hline
			${b_j}$ & 	0& 0& 0& 1& 1& 2& 3& 4& 5& 7& 9& 11& 14& 17& 21\\
			\hline
			$a_j$ &0&0& 0& 1& 2& 4& 7& 11& 16& 23& 32& 43& 57& 74&95
		\end{tabular}
		\bigskip }
	\caption{First values of the sequences $(b_j)$ and $(a_j)$}
	\label{tab:two}
\end{table}

A brief summary on rigid commutators and saturated subgroups 
is presented in  Sec.~\ref{sec:prel}, which also contains a representation of the Sylow $2$-subgroup of the symmetric group
 and to the precise 
definition of the normalizer chain under investigation. The proof of Theorem~\ref{prop:elementsofthen-1normalizer} and 
considerations on unrefinable partitions and rigid commutators can be found in Sec.~\ref{sec:unref}.
 
\section{Preliminaries}\label{sec:prel}
Let $n$ be a non-negative integer.  Let us define a set of permutations $\{s_i \mid 1 \leq i \leq n\}$ which generates a  Sylow $2$-subgroup   of  the  symmetric  group on  $2^n$  letters.
\subsection*{The Sylow $2$-subgroup}
\noindent Let us consider the set
\begin{equation*}
  \mathcal{T}_n=\bigl\{w_1\dots  w_{n}   \mid  w_i  \in
  \{0,1\} \bigr\}
\end{equation*}
of binary words of length $n$, where $\mathcal{T}_0$ contains only
the empty word. The  infinite rooted binary tree  $\mathcal{T}$ is defined
as the graph  whose vertices are $\bigcup_{j\ge  0} \mathcal{T}_j$ and
where two  vertices, say  $w_1\dots w_{n}$  and $v_1\dots  v_{m}$, are
connected   by    an   edge    if   $|m-n|=1$   and    $w_i=v_i$   for
$1 \leq i \leq  \min(m,n)$. The empty word is the  root of the tree
and it is connected with both the two words of length $1$.

\noindent  We  can  define  a   sequence  $\Set{s_i}_{i  \geq  1}$  of
automorphisms of  this tree.  Each  $s_i$ necessarily fixes  the root,
which is the only vertex of degree $2$. The automorphism $s_1$ changes
the  value $w_1$  of the  first letter  of every  non-empty word  into
$\overline{w}_1\deq  (w_1+1)   \bmod  2$  and  leaves   the  other  letters
unchanged.  If $i\ge 2$, we define
\begin{equation*}%\label{eq:generators}
  (w_1\dots w_{k})s_i\deq
  \begin{cases}
    \text{empty word} & \text{if $n=0$} \\
    w_1\dots \overline{w}_i\dots  w_{k} & \text{if $k\ge i$ and $w_1=\dots=w_{i-1}=0$}\\
    w_1\dots w_{k} & \text{otherwise.}
  \end{cases}
\end{equation*}
In general, $s_i$  leaves a word unchanged unless the  word has length
at least $i$ and the letters preceding the $i$-th one are all zero, in
which  case the  $i$-th letter  is increased  by $1$  modulo $2$.   If
$i \le n$  and the word $w_1\dots w_n\in  \mathcal{T}_n$ is identified
with                            the                            integer
$1+\sum_{i=1}^{n}2^{n-i} w_{i}\in \Set{1,\dots, 2^n}$, then $s_i$ acts
on $ \mathcal{T}_n$ as the  the permutation whose cyclic decomposition
is
\begin{equation*}
  \prod_{j=1}^{2^{n-i}}(j,j+2^{n-i})
\end{equation*}
which has order $2$.   In particular, the group $\Span{s_1,\dots,s_n}$
acts  faithfully on the  set  $\mathcal{T}_n$,  whose cardinality  is
$2^n$,  as a  Sylow  $2$-subgroup $\Sigma_n$  of  the symmetric  group
$\Sym(2^n)$. %(see also Fig.~\ref{fig:tree}).
%It is also well known that
%\begin{equation*}
%  \Sigma_{n}	= \Span{s_n} \wr \Sigma_{n-1}   = \Span{s_n} \wr \dots \wr \Span{s_1} \cong \wr_{i=1}^n C_2 
%\end{equation*}  
%is the iterated wreath product of $n$ copies of the cyclic group $C_2$ of order $2$.

\subsection*{Rigid      commutators}     The
\emph{commutator}  of two  elements  $h$ and  $k$ in  a  group $G$  is
defined as $[h,k]\deq h^{-1}k^{-1}hk=h^{-1}h^k$.
The    \emph{left-normed    commutator}    of   the    $m$    elements
$g_1,\dots,g_m\in  G$  is  the  usual   commutator  if  $m=2$  and  is
recursively defined by
\begin{equation*}
  [g_1,\dots,g_{n-1},g_m]\deq \bigl[[g_1,\dots,g_{m-1}],g_m\bigr]
\end{equation*} if
$m\ge 3$.   
In this paper   we  will  only focus  on   left-normed  commutators   in
$s_1, \ldots,  s_n$, therefore, for the sake of simplicity,  we
write   $[i_1,\dots,i_k]$  to   denote   the  left-normed   commutator
$[s_{i_1},\dots, s_{i_k}]$, when $k\ge 2$. We also write $[i]$ to
denote the element $s_i$ and we set \([\ ]\) to be the identity permutation.

\begin{definition}\label{def:rigid_commutators}
  A  left-normed commutator  $r=[i_1,\dots,i_k] \in \Sigma_n$ is  called \emph{rigid}
     when we have $i_1>i_2>\dots >i_k$ or $r=[\ ]$.  The set
  of  all   the  rigid  commutators   of  $\Sigma_n$  is   denoted  by
  $\mathcal{R}$ and 
  we let $\mathcal{R}^*\deq \mathcal{R}\setminus\Set{[\ ]}$.    
\end{definition}

\begin{definition}
	A subgroup of \(\Sigma_n\) is said to be \emph{saturated} if it is generated by rigid commutators. 
\end{definition}
Let us define a special set $\{t_1,t_2,\dots,t_n\}$ of rigid commutators where 
\begin{equation}\label{def:ti}
t_i \deq [i, i-1,\ldots,2,1].
\end{equation}

\begin{remark}
	The saturated subgroup $T\deq\Span{t_1,t_2,\dots,t_n}$
	%= \Span{\puncture{1}{\emptyset},\puncture{2}{\emptyset},\dots,\puncture{n}{\emptyset}}\] 
	is an elementary abelian regular subgroup of \(\Sigma_n\). 
\end{remark}
\begin{theorem}[\cite{aragona2021rigid}.]\label{thm:normalizer_saturated}
	The  normalizer  $N$  in  $\Sigma_n$  of  a  saturated  subgroup  of
	$\Sigma_n$ is also saturated, provided that $N$ contains $T$.
\end{theorem} 

Let $i_1>i_2>\dots >i_k$ and let  $X=\Set{1,\dots,i_1}\setminus \Set{i_1,\dots,i_k}$. For the purposes of this work it is more convenient 
to use the notation
  \[  \puncture{i_1}{X} \] to denote the rigid commutator $[i_1,\dots,i_k]$.
The commutator of two rigid commutators is again a rigid commutator which, in the above notation,  can be computed in the following way:
\begin{proposition}[\cite{aragona2021rigid}]\label{prop:punctured}
	Let  $1  \leq   a,b  \leq  n$  and  let  $I$   and  $J$ be  subsets  of
	$\Set{1,2,\dots,a-1}$ and $\Set{1,2,\dots,b-1}$ respectively. Then
	\begin{equation*}
		\bigl[\,
		\puncture{a}{I} ,
		\puncture{b}{J}
		\, \bigr] = 
		\begin{cases}
			\puncture{\max (a,b)}{\;  (I \cup  J)\setminus\Set{\min(a,b)}} &
			\text{if $\min (a,b)\in I \cup J$}
			\\
			1 & \text{otherwise}.
		\end{cases}
	\end{equation*}
\end{proposition}

\subsection*{The normalizer chain}
The \emph{normalizer chain starting at $T$} is defined as
\begin{equation}\label{eq:normalizer_chain}
  N_n^i \deq \begin{cases}
   N_{\Sigma_n}(T)& \text{if $i=0$},\\
    N_{\Sigma_n}(N^{i-1}_{n}) & \text{if $i\ge 1$.}
  \end{cases}
\end{equation}
Notice that $N_{n}^0$ is the Sylow $2$-subgroup $U_n$ of the affine group $\AGL(2,n)$.
It is worth noticing here that $N_{\Sigma_n}(N_{n}^i)=N_{\Sym(2^n)}(N_{n}^i)$, for all  $i\ge 0$~\cite{Aragona2020}.\\

In order to show where partitions of integers come into play, let us briefly describe the generators of the first $n-2$ normalizers of the chain. First, let
us determine the permutations in  $\Sigma_n$ normalizing $T$:
for      $1      \leq      j<i\leq     n$      let      us      define
\begin{equation}\label{def:uij}
X_{ij}\deq \Set{1,\dots,i}\setminus \Set{j} \text{ and }
  u_{ij}\deq [X_{ij}] = \puncture{i}{\{j\}} \in \mathcal R^*,
\end{equation}
and let us set
\begin{equation*}
  \mathcal{U}_n \deq \Set{t_1,\dots,t_n, u_{ij} \mid 1 \leq j < i \leq
    n} \subseteq \mathcal{R}^*.
\end{equation*}

\noindent Next, let us define
\begin{equation}
  \label{eq:C_lk}
  \mathcal{W}_{ij}\deq \Set{\puncture{i}{I} \in \mathcal{R}^* \,\,\Big  | \,\, I \subseteq
    \Set{1,2,\dots,i-1},  \Size{I}  \ge  2,  \vphantom{\sum}\smash{\sum_{x  \in  I}  x}=j  \;
  }
\end{equation}
for each $1 \leq i \leq n$ and $j$, and
\begin{equation}\label{def:Ni}
  \mathcal{N}_n^{i}\deq
  \begin{cases}
    \mathcal{U}_n & \text{if $i=0$} \\
    \mathcal{N}_n^{i-1}    \dot\cup    \left(    \dot\bigcup_{j=1}^{i}
      \mathcal{W}_{n+j-i,\,j+2} \right) & \text{for $i>0$.}
  \end{cases}
\end{equation}
Note  that, if  $j  \leq i-2$,  then $\Size{  \mathcal{W}_{i,j}}=b_j$,
i.e.\  it corresponds to the number  of partitions  of $j$  into at  least two  distinct
parts. 

The previous elements generate the subgroups in the normalizer chain:
\begin{theorem}[\cite{aragona2021rigid}]\label{thm:normalizers_indices}
  For $i \leq  n-2$, the group $\Span{\mathcal{N}_n^i}$  is the $i$-th
  term $N_n^{i}$ of the normalizer chain. In particular, the subgroup \(N_n^i\) of \(\Sigma_n\) is generated by \(U_n\) and the rigid 
  commutators \(\puncture{a}{X}\) such that 
  \begin{itemize}
  \item \(\Size{X}\ge 2\),
 \item \(\sum _{x\in X}x \le i+2-(n-a) \).
  \end{itemize}
%  i.e.\ all the rigid commutators \(\puncture{a}{X}\) where \(X\) is a partition into distinct parts of a positive integer less or equal to \(i+2-(n-a)\).

\end{theorem}

 The following straightforward consequence is derived:

\begin{corollary}[\cite{aragona2021rigid}]
  For $1\le i\le n-2$, the number $\log_{2}\Size{N^{i}_n : N^{i-1}_n}$
  is  independent  of $n$.   It  equals  the  $(i+2)$-th term  of  the
  sequence $(a_j)$ of the partial sums of the sequence $(b_j)$
  counting the number of partitions of  $j$ into at least two distinct
  parts.
\end{corollary}

\section{Unrefinable partitions}\label{sec:unref}
Every non-empty finite subset $X\subseteq \N\setminus\Set{0}$ represents a
partition of  the integer \[\sum X\deq  \sum_{x\in X} x\]  into \emph{distinct}  parts. Some partitions, e.g.\
$7=1+2+4$, are not refinable, some others are. For example, in $10=1+4+5$, 
the part $5$ can be replaced by $2+3$ obtaining a partition of $10=1+2+3+4$ into distinct parts.
More precisely:
\begin{definition}\label{def:unref}
A finite   non-empty  subset  $X\subseteq  \N\setminus\Set{0}$
 is  \emph{refinable} if there exists $x\in X$
and a  subset  $Y\subseteq \N\setminus  (X \cup  \Set{0})$ such
that              $x=\sum               Y$.         If so, $X$ induces a \emph{refinable partition} into distinct parts of $\sum X$.       Equivalently,
$(X\setminus \Set{x}) \mathbin{\dot\cup} Y$ is again a partition
into distinct  parts of $\sum  X$, called a \emph{refinement}  of $X$.
We  say  that  $X$
is an \emph{unrefinable} partition  of $\sum X$ into distinct parts if $X$ is not refinable.
\end{definition}
\begin{remark}
Notice that if $X$ is a refinable partition,  then there exists $Y$ as in Definition~\ref{def:unref}, $	\Size{Y} = 2$, such that $x=\sum Y$.
\end{remark}
 Let us now recall the concept of minimal excludant, better known from its use in combinatorial game theory~\cite{sprague1935mathematische,grundy1939mathematics,fraenkel2015harnessing} and recently considered in the theory of partitions of integers~\cite{andrews2019,Andrews2020,Ballantine2020}.
\begin{definition}
The \emph{minimal excludant} \(\mex(X)\) of a set $X$ of positive integers is  the least positive integer that is not an element of \(X\).
\end{definition}
\bigskip
Rigid commutators and unrefinable partitions are linked by the following consideration depending on  Proposition~\ref{prop:punctured}.
Suppose  that $a > 1$ and that $X\subseteq  \Set{1,\dots,a-1}$ is refinable. If $x \in X$ and $Y$ are as in Definition~\ref{def:unref},
then 
\[
  [ \puncture{a}{X}, \puncture{x}{Y}]= \puncture{a}{Z}
\]
where $Z \deq (X\setminus \Set{x}) \mathbin{\dot\cup} Y$ is the refinement of $X$ obtained by \emph{replacing} $x$ with $Y$. Notice that  $\sum  Z=\sum X$.   Since  $\sum  Y=x$, by Theorem~\ref{thm:normalizers_indices} the  rigid  commutator
$\puncture{x}{Y}$ belongs to  $N_n^{n-2}$  for   $n\ge   x$.   Conversely,  if
\[
[\puncture{a}{X},     \puncture{x}{Y'}]=     \puncture{a}{Z}\]   for some $Y'$  and if
$\sum       Z=\sum       X$,       then       $x\in       X$       and
$Z = (X\setminus \Set{x}) \mathbin{\dot\cup} Y$  is a  refinement of
$X$, where $Y\deq Y'\setminus X$ is such that $\sum Y=x$.

Bearing this in mind, let us prove the following:
\begin{theorem}\label{prop:elementsofthen-1normalizer}
  Let     $a > 1$, $X \subseteq\{1,2, \ldots,a-1\}$ and    let 
  \(k=\mex(X)\) be the minimal excludant of \(X\).   The rigid
  commutator $\puncture{a}{X}$ belong to  $N_n^{n-1}\setminus N_n^{n-2}$ if and only
  if
  \begin{enumerate}
  \item \label{item:one:refined} $\sum X=a+1$,
  \item \label{item:two:refined} $X$ is an unrefinable partition,
  \item \label{item:three:refined} $a\le n < a+k$.
  \end{enumerate}

\end{theorem}
\begin{proof}
  Suppose that  $c \deq \puncture{a}{X}$ belongs to $N_n^{n-1}\setminus  N_n^{n-2}$ and let
  $b=\sum  X$ and  $h=\min(X)$.   By Theorem~\ref{thm:normalizers_indices}, since $c\notin N_n^{n-2}$, we have $b \geq a+1$.  
  Let us define
  \begin{equation*}
    d \deq \begin{cases}
      u_{h,h-1} & \text{if $h>1$}\\
      t_1 & \text{if $h=1$},
    \end{cases}
  \end{equation*}
  where $u_{h,h-1}$ and $t_1$ are respectively defined in Eq.~\eqref{def:uij} and Eq.~\eqref{def:ti}. 
From Proposition~\ref{prop:punctured}, the  commutator $[c,d]$ is  rigid of  the form  $\puncture{a}{Y}$,
  where $\sum Y=\sum X-1=b-1$.  Since  $d\in N_{n}^{n-2}$, we have that
  $[c,d  ] \in  N_{n}^{n-2}$, and in particular  $\sum Y = b-1\le  a$.   We
  already know that $a \le b-1$, and so we have the equality $\sum X=b=a+1$, as claimed in \eqref{item:one:refined}.

Suppose  now that  $X$ is  refinable. We have already observed  that  there      exist     $x\in      X$     and
  a  rigid commutator $\puncture{x}{Y}\in          N_n^{n-2}$, i.e.\ $\sum Y = x$,          such          that
  \[
  [c,\puncture{x}{Y}] = \puncture{a}{Z},\]
with  $\sum Z=\sum X = a+1$. This
  implies that  $\puncture{a}{Z}\notin N_n^{n-2}$, and
   so $c$  does not
  normalize $N_n^{n-2}$, in contradiction to the assumption. This proves \eqref{item:two:refined}.

  Finally, we obviously have $k\le a$. Let us prove that the equality cannot hold. Suppose indeed that $k=a$, then $X$ is the set $\Set{1,\dots,a-1}$. By~\eqref{item:one:refined},      we      have the quadratic equation
  \[a+1=\sum X= (a^2-a)/2\] which admits no integer solutions.
        Now, let       $l\deq a+k$ and let us prove that $n < l$.    
        Suppose that $n\ge    l$. By definition of $l$ and by Theorem~\ref{thm:normalizers_indices}  we have that
  $d \deq \puncture{l}{\Set{a,k}}\in N_{n}^{n-2}\setminus N_{n}^{n-3}$.      Since    $c\in     N_n^{n-1}$,    then
  $N_{n}^{n-2}\ni  [c,d]= \puncture{l}{X\dot\cup  \Set{k}}$.  This  is
  possible only if $l\ge \sum X  +k =a+1+k =l+1 $, a contradiction. The claim \eqref{item:three:refined} is then proved.\\

  Conversely,      suppose      that    
  \eqref{item:one:refined},        \eqref{item:two:refined}        and
  \eqref{item:three:refined} hold and let us consider a generic rigid commutator in $N_{n}^{n-2}$
  \begin{equation*}
    d \deq\puncture{m}{V}\in N_{n}^{n-2}.
  \end{equation*}
  We have $\sum V \le m$ and $m \le n <a+k$. Let us prove first that $[c,d]\in N_n^{n-2}$. To do so we may assume $[c,d] \neq 1$, otherwise there is nothing to prove. Three cases need to be distinguished.
 
  \noindent Assume   first   that    $m   <   a$.    Since    $[c,d]\ne   1$,   then
  $[c,d]=\puncture{a}{Z}$,  where $Z=(X\setminus\Set{m})\cup  V$.  Notice
  that $\sum Z < \sum X=a+1$, where the  equality is excluded since  $Z$ would otherwise
   represent a refinement  of $X$, which is unrefinable  by hypothesis.  Hence
  $\sum Z \le a$, and so $[c,d]\in N_n^{n-2}$, as claimed.
The case $m=a$ has already been discussed since it corresponds to $[c,d]=1$.
To conclude, suppose  that $a < m \le n <a+k$.   Since $[c,d]  \ne 1$, by Proposition~\ref{prop:punctured} we have $a\in
  V$, and consequently
  \begin{equation}\label{eq:mex}
    \sum (V\setminus \Set{a}) \le m -a < k.
  \end{equation}
 By Eq.~\eqref{eq:mex}, we have that $V\setminus  \Set{a}$ contains only elements smaller than $k$
  and, since $k=\mex(X)$ is  minimum positive integer in  the complementary set
  of $X$ in  $\N$, the set $V\setminus  \Set{a}$ has to be  a subset of
  $X$.      Thus     we     have     $[c,d]=\puncture{m}{Z}$,     where
  $Z=X\cup (V\setminus\Set{a})=X$.  Hence $\sum Z=\sum X=a+1 \le m$ and so $[c,d]\in N_n^{n-2}$.

 We     proved     that      $[c,d]\in     N_n^{n-2}$     for     all
  $d\in N_n^{n-2}\cap \mathcal{R}^*$. Now, since $N_n^{n-2}$ is a saturated
  subgroup       of        $\Sigma_n$      we        have, by Theorem~\ref{thm:normalizer_saturated},   that
  $c\in        N_{\Sigma_n}(N_n^{n-2})=N_{n}^{n-1}$.         Moreover,
  $c\notin N_{n}^{n-2}$ since $\sum X =a+1 >a$.
\end{proof}

Using the \textsf{GAP} implementation of the algorithmic version of Theorem~\ref{thm:normalizers_indices}, available at \textsf{GITHUB} (\url{https://github.com/ngunivaq/normalizer-chain}), it
is easy to compute the first normalizers in the chain of Eq.\eqref{eq:normalizer_chain}. The first values of $\log_2\lvert N_n^{i} : N_n^{i-1}\rvert$ are shown in Table~\ref{tab:4}. Notice that for
		$i\le n-2$ the blue numbers correspond to those of the sequence $(a_j)$, whereas black bold numbers represent $\log_2\lvert N_n^{n-1} : N_n^{n-2}\rvert$ and correspond to the number of unrefinable partitions as in Theorem~\ref{prop:elementsofthen-1normalizer}.
\begin{table}[phbt]
	\centering
	\label{tab:nindices}
	\begin{tabular}{c||c|c|c|c|c|c|c|c|c|c|c|c|c|c}
		$n$ & \multicolumn{14}{c}{  $\vphantom{\Big|} \log_2\lvert N_n^i : N_n^{i-1}\rvert$ for $1 \leq i \leq 14$}                                                \\ \hline \hline 
		3  &  \hl{1} & \hll{0}& 0& 0& 0& 0& 0& 0& 0& 0& 0& 0& 0& 0   \\ \hline 
		4 &  \hl{1} & \hl{2}& \hll{1}& 1& 0& 0& 0& 0& 0& 0& 0& 0& 0& 0   \\ \hline 
		5 &    \hl{1}& \hl{2}& \hl{4}& \hll{1}& 2& 2& 1& 1& 1& 1& 0& 0& 0& 0  \\ \hline 
		6 &   \hl{1}& \hl{2}& \hl{4}& \hl{7}& \hll{2}& 4& 4& 1& 1& 2& 2& 2& 2& 1  \\ \hline 
		7  &   \hl{1}& \hl{2}& \hl{4}& \hl{7}& \hl{11}& \hll{4}& 7& 3& 4& 2& 2& 4& 4& 4  \\ \hline 
		8  &  \hl{1}& \hl{2}& \hl{4}& \hl{7}& \hl{11}& \hl{16}& \hll{7}& 5& 6& 2& 6& 6& 3& 3   \\ \hline 
		9  &   \hl{1}& \hl{2}& \hl{4}& \hl{7}& \hl{11}& \hl{16}& \hl{23}& \hll{4}& 9& 4& 11& 4& 12& 9   \\ \hline 
		10  &   \hl{1}& \hl{2}& \hl{4}& \hl{7}& \hl{11}& \hl{16}& \hl{23}& \hl{32}& \hll{4}& 14& 5& 20& 7& 19   \\ \hline 
		11 &  \hl{1}& \hl{2}& \hl{4}& \hl{7}& \hl{11}& \hl{16}& \hl{23}& \hl{32}& \hl{43}& \hll{5}& 22& 7& 32& 4   \\ \hline 
		12 &   \hl{1}& \hl{2}& \hl{4}& \hl{7}& \hl{11}& \hl{16}& \hl{23}& \hl{32}& \hl{43}& \hl{57}& \hll{7}& 32& 12& 43   \\ \hline 
		13 &   \hl{1}& \hl{2}& \hl{4}& \hl{7}& \hl{11}& \hl{16}& \hl{23}& \hl{32}& \hl{43}& \hl{57}& \hl{74}& \hll{12}& 42& 18   \\ \hline 
		14 &  \hl{1}& \hl{2}& \hl{4}& \hl{7}& \hl{11}& \hl{16}& \hl{23}& \hl{32}& \hl{43}& \hl{57}& \hl{74}&\hl{95}& \hll{8}& 24   \\ \hline 
		15 &   \hl{1}& \hl{2}& \hl{4}& \hl{7}& \hl{11}& \hl{16}& \hl{23}& \hl{32}& \hl{43}& \hl{57}& \hl{74}& \hl{95}& \hl{121}& \hll{8}   \\ 
	\end{tabular}
	\bigskip
	\caption{                   Values                   of
		$\log_2\lvert N_n^i : N_n^{i-1}\rvert$ for  small $i$ and $n$}
	\label{tab:4}
\end{table}
Moreover, notice that the diagonals in the blue area of Table~\ref{tab:4} contain the same sequence of numbers. This is not the case when 
looking at the bold diagonal related to the $(n-1)$-th normalizer. The reason for this is explained below.
\begin{remark}
Recall that, if  \(\sum X \le a \), then the rigid commutator $c = \puncture{a}{X}$ belongs to $N_n^{n-2}$.
In particular $c \in N_m^{m-2}$ for all $m \geq n$. A similar property does not hold for rigid commutators in $N_n^{n-1}$.
  In this case, we can show that if $c \in N_n^{n-1}\setminus N_n^{n-2}$, 
  then $c \in N_m^{m-1}$  for finitely many $m$.
  Indeed, if $a \le  m$, by
  Theorem~\ref{prop:elementsofthen-1normalizer} we have that $c\in N_m^{m-1}$  only if  $m$ is  in the
  interval $a \le  m <a+\mex(X)$. \\
  It is  natural to ask when  $k=\mex(X)$ is the largest
  possible with respect  to $a$. Since $\sum X=a+1$, it  is clear that
  this   happens  exactly   when   $X$  is a triangular partition,   i.e.
  $X=\Set{1,\dots      ,k-1}$.       We       then      have      that
  \[a=\sum X-1=  k(k-1)/2-1=(k+1)(k-2)/2\] is the integer  preceding the
  $(k-1)$-th      triangular       number      $k(k-1)/2$.       Since
  $a=(k+1)(k-2)/2>(k-2)^2/2$ we  have that  $\sqrt{2a}+2$ is  an upper
  bound  for the  largest possible  value  of $k$.  In this case,  if
  $a\ge       3$,       then      $c\in       N_n^{n-1}\setminus N_n^{n-2}$       implies
  \[n    <   a+k    \le    a+\sqrt{2a}+2\le   (\sqrt{a}+1)^2,\]    i.e.\
  $a > n-1-2\sqrt{n}$.
\end{remark}

\begin{example}
When $n=8$ we have computed that $\log_2\lvert N_8^{7} : N_n^{6}\rvert = 7$. Indeed, applying Theorem~\ref{prop:elementsofthen-1normalizer} it can be shown  that the normalizer $N_8^7$ can be generated by the generators of $N_8^6$ and by the following 7 rigid commutators of the type $\puncture{a}{X}$
\begin{eqnarray*}
 \puncture{5}{\{3,2,1\}} &=& [5,4], \\
 \puncture{6}{\{ 4,2,1  \}}&=&  [6,5,3], \\
\puncture{7}{ \{5,2,1\}} &=& [7,6,4,3],\\
\puncture{7}{\{ 4,3,1 \}}&=& [7,6,5,2], \\
 \puncture{8}{\{ 6,2,1  \}}&=&[8,7,5,4,3],\\
\puncture{8}{\{ 5,3,1  \}}&=&[8,7,6,4,2], \\
\puncture{8}{\{  4,3,2 \}}&=&[8,7,6,5,1]. 
\end{eqnarray*}
\end{example}
It is natural to wonder whether a general closed formula for $\log_2\lvert N_n^{n-1} : N_n^{n-2}\rvert$ may be found, even though
at the time of writing this does not seem an easy task.
To our knowledge, the problem of determining a closed formula or a generating function for the sequence $c_j$ of the number of  unrefinable partitions into distinct parts
is open. The  first values of $c_j$ are shown here in Table~\ref{tab:cj}, while a list of the first 1000 integers can be found in the On-Line Encyclopedia of Integer Sequences~\cite[\url{https://oeis.org/A179009}]{OEIS}.
The related problem of determining a closed formula for the number of unrefinable partitions with a given minimal excludant, related to the number of partitions of Theorem~\ref{prop:elementsofthen-1normalizer}, seems at the moment out of reach.

\begin{table}[]
		\centering
	{\renewcommand{\arraystretch}{1.3}
		\begin{tabular}{c||c|c|c|c|c|c|c|c|c|c|c|c|c|c|c}
			$j$ & $0$ & $1$ & $2$ & $3$ & $4$ & $5$ & $6$ & $7$ & $8$ & 9& 10 &11&12&13&14\\
			\hline\hline
			${c_j}$ & 	1& 1& 1& 1& 1& 2& 1& 1& 2& 3& 1& 2& 2& 3&5\\
		\end{tabular}
		\bigskip }
	\caption{First values of the sequences $(c_j)$}
	\label{tab:cj}
\end{table}

\bibliographystyle{amsalpha}
\bibliography{sym2n_ref}

\end{document}